\documentclass [12pt,reqno]{amsart}
\usepackage [normalem]{ulem}
\usepackage {amsmath}
\usepackage {amssymb}
\usepackage {mathrsfs}
\usepackage {enumerate}

\allowdisplaybreaks [2]

\usepackage {color}

\newcommand \Iso {\operatorname {Iso}}

\newcommand \go {\mathcal {G}^{(0)}}

\newcommand {\inv }{^{-1}}

\newcommand {\Z }{\mathbb {Z}}
\newcommand {\N }{\mathbb {N}}

\newcommand {\intiso }{\operatorname {Iso}(\mathcal {G})^\circ }
\newcommand {\group }{\breve {G} \rtimes _{{\breve \rho }} \mathbb {Z}}
\newcommand {\g }{{\bf g}}  

\newcommand {\C }{\mathbb {C}}
\DeclareMathOperator {\supp }{supp}

\newcommand {\B }{{\mathcal {B}}}

\newcommand {\Labe }{{\mathcal {L}}}
\newcommand {\Breg }{{\mathcal {B}_{reg}}}
\newcommand {\Ra }{{\mathcal {R}_\alpha }}
\newcommand {\Ri }{{\mathcal {R}}}
\newcommand {\Clab }{C^*(\B ,\Labe ,\theta )}
\newcommand {\Id }{{\mathcal {I}}}
\newcommand {\He }{{\mathcal {H}}}
\newcommand {\Sa }{{\mathcal {S}}}

\newcommand {\OGE }{{\mathcal {O}_{G,E}}}

\newcommand {\SGE }{{\mathcal {S}_{G,E}}}

\newcommand {\CG }{{\mathcal {G}_{(G,E)}}}

\newtheorem {lemma}[equation]{Lemma}
\newtheorem {cor}[equation]{Corollary}
\newtheorem {thm}[equation]{Theorem}
\newtheorem {prop}[equation]{Proposition}
\newtheorem {noname}[equation]{}

\theoremstyle {definition}
\newtheorem {definition}[equation]{Definition}
\newtheorem {example}[equation]{Example}

\theoremstyle {remark}

\newtheorem *{notation}{Notation}

\numberwithin {equation}{section}

\def\lisa{\textcolor{cyan}}

\begin {document}

\title [Generalised uniqueness theorem for Steinberg algebras]{A Generalised uniqueness theorem and the graded ideal structure of Steinberg algebras}

\author {Lisa Orloff Clark}
\address {Department of Mathematics and Statistics, University of Otago, PO Box 56, Dunedin 9054, New Zealand}
\email {lclark@maths.otago.ac.nz}

\author {Ruy Exel}
\address {Departamento de Matem\'atica, Universidade Federal de Santa Catarina, 88040-970 Florian\'opolis SC, Brazil}
\email {exel@mtm.ufsc.br}\urladdr {http://www.mtm.ufsc.br/~exel/}

\author {Enrique Pardo}
\address {Departamento de Matem\'aticas, Facultad de Ciencias\\ Universidad de C\'adiz, Campus de
Puerto Real\\ 11510 Puerto Real (C\'adiz)\\ Spain.}
\email {enrique.pardo@uca.es}\urladdr {https://sites.google.com/a/gm.uca.es/enrique-pardo-s-home-page/}


\thanks {The first named author was partially supported by Marsden grant 15-UOO-071 from the Royal Society of New Zealand.
The second-named author was partially supported by CNPq. The third-named author was partially supported by PAI III grant FQM-298 of the
Junta de Andaluc\'{\i }a, and by the DGI-MINECO and European Regional Development Fund, jointly, through grant MTM2014-53644-P}




\begin {abstract} Given an ample, Hausdorff groupoid $\mathcal{G}$, and a unital commutative ring $R$, we consider the
Steinberg algebra $A_R(\mathcal {G})$.  First we prove a uniqueness theorem for this algebra and then, when $\mathcal{G}$ is graded by a cocycle, we study
graded ideals in $A_R(\mathcal {G})$.    Applications are given for two classes of ample groupoids, namely those coming from
actions of groups on graphs, and also to groupoids  defined in terms of Boolean dynamical systems.
\end {abstract}

\maketitle

\section {Introduction}
Steinberg algebras were independently introduced in \cite {Steinberg} and \cite {CFST} and are closely related
to the machinery established in \cite {E}. This broad class of algebras
provides a general model for Leavitt path algebras associated to directed graphs, Kumjian Pask algebras
associated to higher-rank graphs, and discrete inverse semigroup
algebras.

In this note, we use an analysis of the `interior of the isotropy bundle' to prove a generalised uniqueness theorem.  We
also characterize the graded ideals for Steinberg algebras over groupoids equipped with a cocycle taking values in a
discrete group.  Applications are then provided in two broad classes of examples.  Our presentation is as follows:

After establishing our notation and discussing some preliminary results in Section \ref {sec:prelim}, we move to
Section \ref {sec:gut} where we prove a generalised uniqueness theorem for Steinberg algebras modeled after the
$C^{\ast }$-algebra \cite [Theorem~3.1(c)]{BNRSW}.  Our uniqueness theorem says that a Steinberg algebra homomorphism is
injective if and only if it is injective on the interior of the isotropy group bundle.  This generalises theorems \cite
[Theorem~5.2]{GN} and \cite [Theorem~5.4]{CGN} for Leavitt path algebras and Kumjian-Pask algebras respectively.

In Section~\ref {sec:gi}, we characterize the graded ideals of a Steinberg algebra built from `graded groupoids'.
These are groupoids that come equipped with a homomorphism (or `cocycle') $c$ into a discrete group such that the inverse image of the
identity doesn't have too much isotropy.  In this setting, we show that the graded ideals are precisely those generated by
open invariant subsets of the unit space.

In the last two sections, we give two more classes of examples.  Both are classes whose $C^{\ast }$-algebras have been
constructed using the machinery of \cite {E}.  In Section~\ref {sec:gag} we consider groups acting on graphs, as defined in \cite {EP_GGS}; we
describe the Steinberg algebra associated associated to such an action. This broad class of algebras includes the class of Leavitt path
algebras.   We also get an algebraic analogue of Katsura's algebras $\mathcal {O}_{A,B}$.
Finally, in Section~\ref {sec:bds} we construct a Steinberg algebra from a Boolean dynamical system, as introduced in \cite {COP}.



  \newdimen \boxht
  \def \sumTwoIndices #1#2{\sum _{\buildrel {\scriptstyle #1}\over {\vrule height 8pt width 0pt #2}}}
  \def \med #1{\mathop {\textstyle #1}\limits }
  \def \medoplus {\med \bigoplus }

\section {Preliminaries}
\label {sec:prelim}

Recall that a
  topological\footnote {Unless otherwise mentioned, all topological spaces (including topological groupoids) in this work
 \lisa{are} assumed to be Hausdorff.}
  space $X$ is said to be \emph {zero-dimensional} if the topology of $X$ admits a basis
consisting of  clopen (closed and open) sets.  If, in addition, $X$ is locally compact, it is easy to see
that $X$ also admits a basis formed by \emph {compact} open sets.

In this work, we assume most topological spaces are locally compact, zero-dimensional and Hausdorff.  
If $X$ is such a space, and if $R$ is a unital commutative ring, we denote by
  $$
  C_c(X,R)
  $$
  the
  $R$-module\footnote {Even though $C_c(X,R)$ is also an $R$-algebra, relative to pointwise product, we will not always view it
as such.}
  formed by all locally constant, compactly supported, $R$-valued functions on $X$.  Notice that an $R$-valued function is
locally constant if and only if it is continuous once we equip $R$ with the discrete topology.

If $D$ is a compact open subset of $X$, the characteristic function of $D$, here denoted by $1_D$, is clearly an element
of $C_c(X,R)$.  In fact, every $f$ in $C_c(X,R)$ may be written as
  $$
  f=\sum _{i=1}^n a_i1_{D_i}
  $$
  where the $D_i$ are compact open, pairwise disjoint subsets of $X$.  The \emph {support} of $f$, defined by
  $$
  \text {supp}(f) = \{x\in X: f(x)\neq 0\},
  $$
  (we do not use closure in the definition of the support), is clearly a compact open subset.

If  $X$ is as above, and $U\subseteq X$ is an open subset, then $U$ is also  a locally compact, zero-dimensional, Hausdorff space.  Moreover we may view
$C_c(U,R)$ as a submodule of $C_c(X,R)$ as follows:  given any $f$ in $C_c(U,R)$, let $\tilde f$ be the extension  of $f$
to
$X$ defined to be zero on $X\setminus U$.  It is then easy to see that $\tilde f$ is continuous on $X$ and, in fact, that
it lies in $C_c(X,R)$.
  The correspondence
  \begin {equation}
  \label {equation:NatInclusion}
  f\in C_c(U,R) \mapsto \tilde f \in C_c(X,R),
  \end {equation}
  is then an $R$-module isomorphism from $C_c(U,R)$ onto the submodule of $C_c(X,R)$ formed by the functions vanishing on
$X\setminus U$.  We may therefore view $C_c(U,R)$ as a submodule of $C_c(X,R)$.

Given a groupoid ${\mathcal G}$, we will always denote its \emph {unit space} by $\go $, the set of \emph {composable pairs} by ${\mathcal G}^{(2)}$, and its
\emph {source} and \emph {range} maps by $s$ and $r$, respectively.
A \emph {bissection} in ${\mathcal G}$ is a subset $U\subseteq {\mathcal G}$ such that the restrictions of $r$ and $s$ to $U$ are both injective.

A Hausdorff topological groupoid ${\mathcal G}$ is said to be \emph {\'etale} if
  $\go $ is locally compact and Hausdorff in the relative topology, and
  its range map is a local homeomorphism from ${\mathcal G}$ to $\go $ (the source map will consequently share that property).
  It is easy to see that the topology of an \'etale groupoid admits a basis formed by open bissections.
  In an \'etale groupoid one has that $\go $ is open in ${\mathcal G}$.
If, in addition,  ${\mathcal G}$ is Hausdorff, then $\go $ is also closed in ${\mathcal G}$.  Throughout
this paper all groupoids will be assumed Hausdorff.

\begin {definition}
An \'etale,  Hausdorff groupoid ${\mathcal G}$ is said to be \emph {ample} if   $\go $ is zero-dimen\-sional.
\end {definition}

It is easy to see that the topology of
an ample groupoid admits a basis formed by \emph {compact open bissections}.
  En passant we deduce that an ample groupoid
is also  locally compact and zero-dimensional.

\begin {definition}[\cite {Steinberg}, \cite {CFST}]
Given an ample,
  Hausdorff\footnote {Steinberg algebras may also be defined for ample groupoids that are not Hausdorff
\cite [Definition 4.1]{Steinberg}, but not as a straightforward generalization of the above.}
  groupoid ${\mathcal G}$, and a unital commutative ring $R$, the \emph {Steinberg algebra} associated to
${\mathcal G}$, denoted $A_R(\mathcal G)$, is defined to be the $R$-algebra obtained by equipping
$C_c({\mathcal G},R)$ with
the \emph {convolution product}
  $$
  (fg)(\gamma ) = \sumTwoIndices {(\gamma _1,\gamma _2)\in {\mathcal G}^{(2)}}{\gamma _1\gamma _2=\gamma }f(\gamma _1)g(\gamma _2).
  $$
\end {definition}

Since every element of $A_R(\mathcal G)$ is a linear combination of characteristic functions of bissections, it is
interesting to notice that
  \begin {equation}
  \label {equation:ProdCharBissect}
  1_B  1_D = 1_{BD},
  \end {equation}
  whenever $B$ and $D$ are compact open bissections in ${\mathcal G}$.

If ${\mathcal G}$ is an ample, Hausdorff groupoid, and if ${\mathcal H}$ is an open subgroupoid of ${\mathcal G}$,
then ${\mathcal H}$ is also ample and the natural inclusion
  \begin {equation}
  \label {equation:NatAlgInclusion}
  A_R(\mathcal H)\to A_R(\mathcal G),
  \end {equation}
  given by (\ref {equation:NatInclusion}), identifies
  $A_R(\mathcal H)$ with a subalgebra of  $A_R(\mathcal G)$.

An elementary example  of an open subgroupoid of ${\mathcal G}$ is $\go $, so we have that
  $$
  A_R(\go) \subseteq A_R(\mathcal G).
  $$
  Due to the very special nature of the groupoid $\go $, the convolution product on
$A_R(\go) $ coincides with the pointwise product.


Recall that an algebra $B$ is
said to be $s$-\emph {unital} when for every
  $
  r_1, r_2, ..., r_n \in R,
  $
  we can find $s\in R$ such that
  $$
  sr_i = r_i = r_is,\quad \forall i=1,\ldots ,n.
  $$
  Moreover, if $E$ is a given set of idempotent elements in $B$ such that the above $s$ can be taken in $E$, we will say that $E$ is a set
of \emph {local units} for $B$.

\begin {lemma}\label {equation:SteinbSUnital}
  Let ${\mathcal G}$ be a ample, Hausdorff groupoid.  Then the family of all idempotent elements of   $A_R(\go) $ is a set of local
units for
$A_R(\mathcal G)$.   In particular the latter is an $s$-unital algebra.
\end {lemma}

\begin {proof}
  Given a compact open bissection  $D\subseteq {\mathcal G}$, and given any compact open subset $C\subseteq \go $, with $s(D)\cup r(D)\subseteq C$, one has
that   $DC=D = CD$, so
  $$
  1_D1_C=1_D = 1_C1_D.
  $$
  Thus, if we are given elements
 $$
  f_i = \sum _{D \in F_i} a_D1_D\in A_R(\mathcal G), \quad 1\leq i\leq n,
  $$
  where the $F_i$ are finite sets of compact open bissections, we may take
  $$
  C=\bigcup _{1\leq i\leq n}\bigcup _{D\in F_i}s(D)\cup r(D),
  $$
  and it will follow that $1_C$ is an idempotent element of $A_R(\go) $, and
  $$
  f_i1_C=f_i = 1_Cf_i, \quad \forall i=1,\ldots ,n.
  \eqno \qedhere
  $$
\end {proof}

An immediate consequence of the above result is that
  \begin {equation}
  \label {equation:NonGegenSubAlg}
  A_R(\mathcal G)= A_R(\go) A_R(\mathcal G)=   A_R(\mathcal G)A_R(\go) .
  \end {equation}

Recall that a subset $U \subseteq \go $ is said to be \emph {invariant} if for every $\gamma $ in ${\mathcal G}$ one has that
  \begin {equation}
  \label {equation:InvariantSubset}
  r(\gamma )\in U  \iff s(\gamma )\in U.
  \end {equation}
  Given such an invariant subset $U$, let
  \begin {eqnarray*}
  {\mathcal G}|_U
    & = & \{\gamma \in {\mathcal G}: r(\gamma )\in U\} \\
    & = & \{\gamma \in {\mathcal G}: s(\gamma )\in U\},
  \end {eqnarray*}
  and observe that   ${\mathcal G}|_U$ is a subgroupoid of ${\mathcal G}$.  If, in addition,  $U$ is an open subset of $\go $, then
  ${\mathcal G}|_U$ is open in ${\mathcal G}$.   Therefore, as in \eqref {equation:NatAlgInclusion}, we have that
$A_R({\mathcal G}|_U)$ is a subalgebra of  $A_R(\mathcal G)$.

Regardless of whether or not $U$ is invariant, as long as it is an open subset of $\go $ we may always see $U$ itself as
an open subgroupoid of ${\mathcal G}$.  So,
again by \eqref {equation:NatAlgInclusion}, we may identify
$A_R(U)$ as a subalgebra of  $A_R(\mathcal G)$, observing that $A_R(U)$ is nothing but $C_c(R,U)$ with pointwise multiplication.

Given an algebra $A$, and given two subsets $X$ and $Y$ of $A$, we will denote by $XY$ the
$R$-submodule of $A$ generated by the set of products, namely
  $$
  XY = \big \{\sum _{i=1}^n r_ix_iy_i: n\in \mathbb N,\ r_1,\ldots ,r_n\in R,\ x_1,\ldots ,x_n\in X,\ y_1,\ldots ,y_n\in Y\big \}.
  $$

  This notation is used  in the next:

\begin {prop}\label {equation:InvarCriteria}
Let ${\mathcal G}$ be an ample, Hausdorff groupoid. Given an open subset $U\subseteq \go $, one has that $U$ is invariant in the sense of \eqref {equation:InvariantSubset}
if and only if
  $$
  \phantom {A_R({\mathcal G}|_U) = }  A_R(\mathcal G)A_R(U)=  A_R(U)  A_R(\mathcal G).
  $$
  In this case, one has that
  $$
  A_R({\mathcal G}|_U) = A_R(\mathcal G)A_R(U) =  A_R(U)  A_R(\mathcal G),
  $$
  and consequently $A_R({\mathcal G}|_U)$ is an ideal in $A_R(\mathcal G)$.
\end {prop}

\begin {proof}
   Supposing that $U$ is invariant,
   let us first prove that  $A_R(\mathcal G)A_R(U) \subseteq A_R({\mathcal G}|_U)$. For this it suffices to show that, given any
compact open bissection $B\subseteq {\mathcal G}$, and any compact open subset $K\subseteq U$, one has that $1_B1_K\in A_R({\mathcal G}|_U)$.  To see this,
notice that
  $
  1_B1_K = 1_{BK},
  $
  as seen in   \eqref {equation:ProdCharBissect},
  and that if $\gamma $ is any element of ${\mathcal G}$ lying in $BK$, then $s(\gamma )\in K\subseteq U$, whence $\gamma \in {\mathcal G}|_U$.  This shows that $BK\subseteq {\mathcal G}|_U$, so $1_{BK}$
indeed belongs to $A_R({\mathcal G}|_U)$.

   Let us now prove that  $A_R(\mathcal G)A_R(U) \supseteq A_R({\mathcal G}|_U)$.  Again we may focus on characteristic functions, meaning
that all we
must do is show that $1_B$ lies in $A_R(\mathcal G)A_R(U)$, for all compact open bissections $B\subseteq {\mathcal G}|_U$.
Given such a $B$, observe that $s(B)\subseteq U$, so $1_{s(B)}\in A_R(U)$, and
  $$
  1_B =
  1_{Bs(B)} =
  1_B1_{s(B)} \in
  A_R(\mathcal G)A_R(U).
  $$
  This proves that   $A_R({\mathcal G}|_U) = A_R(\mathcal G)A_R(U)$, and one may similarly prove that
  $
  A_R({\mathcal G}|_U) =A_R(U)  A_R(\mathcal G),
  $
  which readily implies that $A_R({\mathcal G}|_U)$ is an ideal in $A_R(\mathcal G)$.

In order to complete the proof we must show that the first identity displayed in the statement is a sufficient condition
for the invariance of $U$.  So, pick any $\gamma \in {\mathcal G}$, such that $s(\gamma )\in U$, and let us prove that  $r(\gamma )\in U$, as well.  For
this, choose a compact open bissection $B$ containing $\gamma $. Replacing $B$ with $BK$, where $K$ is any compact open
set such that
  $$
  s(\gamma ) \in K\subseteq U,
  $$
  we may suppose that $s(B)\subseteq U$.  It follows that
  $$
  1_B = 1_B1_{s(B)} \in A_R(\mathcal G)A_R(U)=  A_R(U)  A_R(\mathcal G).
  $$
  Since $A_R(U)$ has local units by \eqref {equation:SteinbSUnital}, there exists an idempotent element in $A_R(U)$,
necessarily of the form $1_V$, for some compact open set $V\subseteq U$, such that $1_B=1_V1_B$.  It follows that $B=VB$, and
hence
  $$
  r(\gamma ) \in r(B) = r(VB) \subseteq V \subseteq U.
  $$
  Similarly one proves that ``$r(\gamma )\in U\Rightarrow s(\gamma )\in U$", so $U$ is invariant.
\end {proof}

  Given an ample, Hausdorff groupoid ${\mathcal G}$,
 for $u \in \go $, we let ${\mathcal G}^u_u$ denote the isotropy group at $u$.
 The isotropy group bundle is the set
  \[ \Iso ({\mathcal G}) := \bigcup _{u \in \go }{\mathcal G}^u_u.\]
   We write $\intiso $ for the interior of the isotropy group bundle 
in the relative topology.  Since $\go $ is open, one has that
 $\go \subseteq \intiso $.

This gives us two canonical examples of open subgroupoids of ${\mathcal G}$, namely $\go $ and $\intiso $.  In the spirit of
\eqref {equation:NatAlgInclusion} we then have
  $$
  A_R(\go) \subseteq A_R(\intiso) \subseteq A_R(\mathcal G).
  $$

 We say ${\mathcal G}$ is \emph {principal} if $\operatorname {Iso}({\mathcal G})=\go $.
 One can show this is equivalent to saying that the map $\gamma \mapsto (r(\gamma ),s(\gamma ))$ is injective.
We call a groupoid ${\mathcal G}$ \emph {effective}\footnote {
This has also been called \emph {essentially principal} in a number of papers and that term is also somewhat standard.
  However, there are papers where essentially principal is used to  mean something different
  so we have chosen to use `effective' to avoid
  confusion.}
 if $\intiso = \go $.
However, if $C$ is a closed invariant subset of $\go $, the restriction ${\mathcal G}|_C$ might not be effective, even if ${\mathcal G}$
has this property.
  We therefore say that ${\mathcal G}$ is \emph {strongly effective} if ${\mathcal G}|_C$ is effective for every closed invariant subset $C\subseteq \go $.
\section {A General uniqueness theorem}
\label {sec:gut}

In this section, we prove a generalised uniqueness theorem for
the Steinberg algebra $A_R({\mathcal G})$ associated to an ample, Hausdorff groupoid.
Our theorem is the algebraic
analogue of \cite [Theorem~3.1(c)]{BNRSW} for groupoid $C^{\ast }$-algebras.

\begin {thm}[Generalized Uniqueness Theorem]
\label {thm:Uniqueness}
Let ${\mathcal G}$ be a second-countable, ample, Hausdorff groupoid and
let $R$ be a unital commutative ring.
Suppose that $A$ is an $R$-algebra and that $\pi :A_R({\mathcal G}) \to A$ is a ring homomorphism.
Then $\pi $ is injective if and only if $\pi \circ \iota $ is injective.
\end {thm}

 We summarize the results we use to prove Theorem~\ref {thm:Uniqueness}
 in the following two lemmas.

\begin{lemma}[{\cite [Lemma~3.3(a)]{BNRSW}}]
\label{lem:denseuu}
Let ${\mathcal G}$ be a second-countable, ample, Hausdorff groupoid.
Then \[X := \{u \in \go : {\mathcal G}^u_u \subseteq \intiso \}\]
is dense in $\go $.
\end {lemma}

The next lemma is the algebraic analogue of \cite [Lemma~3.3(b)]{BNRSW}.
We check that it holds in a purely algebraic setting.

\begin{lemma}
 \label{lem:getKu}
 Let ${\mathcal G}$ be a second-countable, ample, Hausdorff groupoid
 and
let $R$ be a unital commutative ring.  Suppose $u \in \go $ is such that
 ${\mathcal G}^u_u \subseteq \intiso $ and take $f \in A_R({\mathcal G})$ such that
 there exists $\gamma _u \in {\mathcal G}^u_u$ with
 $f(\gamma _u) \neq 0$.  Then
 there exists a compact open set $K \subseteq \go $ such that $u \in K$ and
 \[0 \neq 1_Kf1_K \in \iota (A_R(\intiso )).\]
\end{lemma}

\begin{proof}
 Write $f= \sum _{D \in F}a_D1_D$ where $D$ is a disjoint collection of compact open bissections.  For each $D \in F$,
 choose a compact open neighbourhood $V_D \subseteq \go $ as follows:
 \begin{itemize}
  \item If $u = r(\gamma ) = s(\gamma )$ for some $\gamma \in D$, then
  $\gamma \in \intiso $ by assumption.  Let
  $V_D$ be a compact open subset in $\go $ containing $u$ such that $V_D$
  is contained in the open set \[r(D \cap \intiso ) = s(D \cap \intiso ).\]
Then  $V_D D V_D \subseteq D \cap \intiso $.
\item If there exists $\gamma \in D$ such that $r(\gamma ) = u$ and $s(\gamma ) \neq u$ or $s(\gamma ) = u$ and
$r(\gamma ) \neq u$, we chose $V_D$ as follows.  Because ${\mathcal G}$ is Hausdorff, we can find a compact open subset
$D' \subseteq D$ containing $\gamma $ such that $r(D') \cap s(D') = \emptyset $.
Take $V_D = r(D')$ (or $V_D= s(D')$), so that $u \in V_D$ and  $V_D DV_D = \emptyset $.
\item If $u \notin r(D)$ and $u \notin s(D)$, use that ${\mathcal G}$ is Hausdorff to choose a neighbourhood
$V_D$ of $u$ such that $V_D DV_D=\emptyset $.
 \end {itemize}
Let \[K:= \bigcap _{D \in F} V_D.\]  Then $K$ is a compact open set that contains $u$ and by construction
\[1_Kf1_K \in \iota (A_R(\intiso )).\]
  Furthermore
  $$
  1_Kf1_K(\gamma _u) = 1_K(r(\gamma _u))f(\gamma _u)1_K(s(\gamma _u)) = f(\gamma _u) \neq 0,
  $$
  so $1_Kf1_K \neq 0$.
\end {proof}

\begin {proof}[Proof of Theorem~\ref {thm:Uniqueness}]
 We show that $\pi $ not injective implies $\pi \circ \iota $ not injective.
 Suppose $\ker \pi \neq \{0\}$.
 By \cite [Lemma~3.1]{CE} there exists a function $f \in \ker \pi $ such that
 $\supp f \cap \go \neq \emptyset $.
 Because ${\mathcal G}$ is Hausdorff,  $\supp f \cap \go $ is an open subset of $\go $.
 Thus Lemma~\ref {lem:denseuu} implies that there exists $u \in \supp f \cap \go $
 such that  ${\mathcal G}^u_u \subseteq \intiso $.
Now apply
Lemma~\ref{lem:getKu} (with $\gamma _u=u$) to get a nonempty compact open set $K \subseteq \go $ such that
\[\supp (1_Kf1_K) \subseteq \iota (A_R(\intiso )).\]
Since $\ker \pi $ is an ideal, we also have $1_K f 1_K \in \ker \pi $.
Thus $0\neq 1_Kf1_K \in \ker (\pi \circ \iota )$.
\end {proof}

In Corollary~\ref {cor:ckuniquness} below, we see that
the
usual Cuntz-Krieger uniqueness theorem (see \cite [Theorem~3.2]{CE} or \cite [Proposition~3.3]{St2016}) for effective
ample groupoids is a consequence of Theorem~\ref {thm:Uniqueness}.
In particular, if ${\mathcal G}$ is effective, then $\intiso = \go $.
Thus  $\pi \circ \iota $ is injective if and only if $\pi (r1_K) \neq 0$ for all compact $K \subseteq \go $ and $r \in R \setminus \{0\}$.

\begin {cor}[The Cuntz-Krieger Uniqueness Theorem]
 \label {cor:ckuniquness}
 Let ${\mathcal G}$ be an effective, second-countable, ample, Hausdorff groupoid and
let $R$ be a unital commutative ring.
Suppose $A$ is an $R$-algebra and $\pi :A_R({\mathcal G}) \to A$ is a ring homomorphism.
Then $\pi $ is injective if and only if $\pi (r1_K) \neq 0$ for all compact open subset $K \subseteq \go $ and $r \in R \setminus \{0\}$.
\end {cor}

 In \cite[Proposition 2.3]{E2011} it is shown that the Cuntz-Krieger uniqueness theorem
 is false for $C^{\ast }$-algebras associated to non-Hausdorff groupoids.
 Similarly, it is false in for the Steinberg algebras associated to ample, non-Hausdorff groupoids.
Although both Lemma~\ref{lem:denseuu} and Lemma~\ref{lem:getKu} hold in the non-Hausdorff setting,
the following slight modification of \cite[Proposition 2.3]{E2011} illustrates that 
the uniqueness theorem itself fails. 

\begin {example}
Let ${\mathcal G}$ be the groupoid described in
\cite [Section 2]{E2011}, and recall that the unit space of ${\mathcal G}$ is the subset $Z\subseteq {\mathbb R}^2$ obtained as the union
of the sets
  $$
  X  = [-1, 1]\times \{0\},  \hbox {\quad and\quad }  Y = \{0\}\times [-1, 1].
  $$

  Let $L$ be the subset of $[-1,1]$ given by
  $$
  L=\{\pm 1/n: n\in {\mathbb N},\ n\geq 1\}  \cup \{0\},
  $$
  and put
  $$
  X'  = L\times \{0\},  \hbox {\quad and\quad }  Y' = \{0\}\times L.
  $$
  Then $Z':=X'\cup Y'$ is easily seen to be a closed invariant subset of $Z$, so
  $$
  {\mathcal G}':={\mathcal G}|_{Z'}
  $$
  is an \'etale subgroupoid of ${\mathcal G}$.  Observing that $Z'$ is zero-dimensional, we have that ${\mathcal G}'$ is an ample groupoid and
it is not hard to see that it is also effective.

If we replace ${\mathcal G}$ by ${\mathcal G}'$, and $f$ by its restriction to ${\mathcal G}'$ in \cite [Proposition 2.2]{E2011}, then it is clear
that the conclusion of this result still holds, and hence, as in \cite [Proposition 2.3]{E2011}, one has that $\mathbb Cf$ is a nontrivial ideal in
$A_{\mathbb C}({\mathcal G}')$ trivially intersecting $A_{\mathbb C}({\mathcal G}'\null ^{(0)})$.
 \end {example}
\section {Ideals in graded algebras}

Our next major goal will be to study graded ideals in
Steinberg algebras.
In preparation for this
we will now discuss certain aspects of ideals in abstract graded algebras.
Throughout this section we will fix a unital commutative ring $R$.  All algebras in this section will be understood to
be $R$-algebras.

\begin {definition}
  We will say that an algebra $A$ is graded over a given group $\Gamma $, or that $A$ is $\Gamma $-graded, if $A$ is equipped with a direct sum decomposition
  $$
  A = \bigoplus _{g\in \Gamma }A_g,
  $$
  where each $A_g$ is an $R$-submodule, satisfying
  $$
  A_gA_h\subseteq A_{gh}, \quad \forall g,h\in \Gamma .
  $$
  In this case we will say each $A_g$ is the \emph {homogeneous
submodule} associated to the group element $g$.
\end {definition}

The central concept to be studied in this section is defined next.

\begin {definition}
  Let $A$ be a $\Gamma $-graded algebra.
  An ideal $J\trianglelefteq A$ is said to be a \emph {graded ideal}\/ if $J = \bigoplus _{g\in \Gamma }J_g$, where each $J_g=J\cap A_g$.
\end {definition}

It is an elementary exercise to show that $J$ is graded if and only if $\pi _g(J)\subseteq J$, for every $g$ in $G$, where $\pi _g$ is
the standard projection of $A$ onto $A_g$.

\begin {definition}\label {def:Invariant}
  Given a $\Gamma $-graded algebra, we will say that an ideal $I\trianglelefteq A_e$ is
\emph {invariant} provided
  $$
  A_gI = IA_g, \quad \forall g\in \Gamma .
  $$
\end {definition}

Invariant ideals give rise to graded ideals as follows:

\begin {prop}
  Given a  $\Gamma $-graded algebra, and an invariant ideal $I\trianglelefteq A_e$, one has that
  $$
  \text {Ind}(I) = \bigoplus _{g\in \Gamma }IA_g,
  $$
  is a graded ideal in $A$, henceforth called the ideal \emph {induced} by $I$.
\end {prop}

\begin {proof}
  Left to the reader.
\end {proof}

It would be nice if invariant ideals in $A_e$ were in a bijective correspondence with graded
ideals in $A$, according to the above procedure, but this is not always true, unless we restrict ourselves to special
situations which we  now describe.

In order to introduce our first major condition, observe that the polynomial algebra $R[X]$ is a
$\mathbb {Z}$-graded algebra; the homogeneous submodule associated to a nonnegative integer $n$ is the set of all
monomials $aX^n$, for $a$ in $R$, while all other homogeneous submodules are zero.  This disparity between
  homogeneous submodules for positive and negative integers is an undesirable feature of this graded
algebra which we will rule out via the following:

\begin {definition}
  We will say that the graded algebra   $A = \bigoplus _{g\in \Gamma }A_g$, is \emph {symmetrically graded} if
  $$
  A_gA_{g\inv }A_g = A_g,
  $$
  for each $g$ in $\Gamma $.
\end {definition}

Observe that a necessary condition for  $A = \bigoplus _{g\in \Gamma }A_g$ to be a symmetrically graded algebra is that
$A_e=A_e^3$, which implies that
  $$
  A_e =  A_e^3 \subseteq A_e^2 \subseteq A_e,
  $$
  so that   $A_e = A_e^2$, which is to say that $A_e$ is an idempotent algebra.

Speaking of graded ideals, observe that these may also be considered as graded algebras which may or may  not be
symmetrically graded.  For the case of induced ideals we have:

\begin {prop}\label {prop:IndIsSymGraded}
  Given a symmetrically\/ $\Gamma $-graded algebra $A$, and an invariant ideal $I\trianglelefteq A_e$, consider
the following statements:
  \begin {enumerate}[(a)]
  \item $I$ is idempotent,
  \item $\text {Ind}(I)$ is symmetrically graded.
  \end {enumerate}
  Then (a) implies (b) and the converse holds provided $A_e$ is $s$-unital.
  \end {prop}

\begin {proof}
  Set $J=\text {Ind}(I)$, so $J_g = J\cap A_g = IA_g$.  Assuming $(a)$, we then have
  $$
  J_gJ_{g\inv }J_g = IA_gIA_{g\inv }IA_g = I^3A_gA_{g\inv }A_g = IA_g = J_g,
  $$
  so $J$ is symmetrically graded.  Conversely, assuming $(b)$, we have already seen that $J_e$ is idempotent.  Moreover,
since $A_e$ is $s$-unital we have
  $$
  I = A_eI = J_e,
  $$
  from where $(a)$ follows.
\end {proof}

The above result may be seen as a process to produce symmetrically graded ideals in $A$ out of invariant idempotent
ideals of $A_e$, and we will now prove that all
such ideals arise from this construction.

\begin {prop}\label {prop:AllSymGrdIdealsAreInd}
  Given a $\Gamma $-graded algebra, and a graded ideal $J\trianglelefteq A$, suppose that $J$ is
symmetrically graded.  Then $J_e=J\cap A_e$ is an invariant idempotent ideal of $A_e$, and $J = \text {Ind}(J_e)$.
\end {prop}

\begin {proof} Setting $J_g=J\cap A_g$,  we have
  $$
  J_eA_g =
  (J\cap A_e)A_g \subseteq
  J\cap A_g =
  J_g =
  J_gJ_{g\inv }J_g \subseteq
  J_eJ_g \subseteq
  J_eA_g,
  $$
  from where we see that $J_eA_g = J_g$ and one may similarly prove that $A_gJ_e=J_g$.  In particular $J_eA_g = A_gJ_e$,
so $J_e$ is invariant and clearly idempotent.  Moreover, since $J$ is graded, we have
  $$
  J=  \bigoplus _{g\in \Gamma }J_g=  \bigoplus _{g\in \Gamma }J_eA_g=\text {Ind}(J_e).
  \eqno {\qedhere }
  $$
  \end {proof}

Aided by our last results we may now prove:

\begin {prop}\label {prop:GradedCorrespondence}
  Let $A = \bigoplus _{g\in \Gamma }A_g$ be a symmetrically graded algebra such that $A_e$ is $s$-unital.  Then the correspondence
  $$
  I \mapsto \text {Ind}(I)
  $$
  is a bijection from the set of all invariant idempotent ideals $I\trianglelefteq A_e$ onto the set of all symmetrically
graded ideals $J\trianglelefteq A$.
\end {prop}

\begin {proof}
Given an invariant idempotent ideal $I\trianglelefteq A_e$, we have by (\ref {prop:IndIsSymGraded}) that $\text {Ind}(I)$ is
indeed a symmetrically graded ideal, hence our correspondence is well defined. It is moreover surjective by
(\ref {prop:AllSymGrdIdealsAreInd}).

Given any invariant ideal $I\trianglelefteq A_e$, observe that
  $$
  \text {Ind}(I)\cap A_e = IA_e = I,
  $$
  where the last equality is a consequence of the fact that $A_e$ is $s$-unital.    If $I'$ is another invariant ideal with   $\text {Ind}(I)= \text {Ind}(I')$, then
  $$
  I = \text {Ind}(I)\cap A_e = \text {Ind}(I')\cap A_e = I'.
  $$
  This shows that $\text {\sl Ind}$ is one-to-one, so our proof is concluded.
\end {proof}
\section {Graded ideals in Steinberg algebras}
\label {sec:gi}

Throughout this section ${\mathcal G}$ will always refer to an ample, Hausdorff groupoid.
  Given a group $\Gamma $,   recall that a map
  $$
  c:{\mathcal G}\to \Gamma ,
  $$
  is said to be a \emph {cocycle} if $c(\gamma _1\gamma _2)=c(\gamma _1)c(\gamma _2)$, for every $(\gamma _1,\gamma _2)\in {\mathcal G}^{(2)}$.  Cocycles are sometimes also called homomorphisms.
Given such a cocycle, for each $g\in \Gamma $, we shall denote by
  $$
  {\mathcal G}_g = \{\gamma \in {\mathcal G}: c(\gamma )=g\}.
  $$
  Assuming that $c$ is locally constant (which is to say that $c$ is continuous when $\Gamma $ is given the discrete
topology), observe that each ${\mathcal G}_g$ is an open set.  We may then view $C_c({\mathcal G}_g,R)$ as a submodule of $A_R(\mathcal G)$, as in
\eqref {equation:NatInclusion}.

\begin {prop}\label {equation:SteinSymmGrded}
  Given an  ample, Hausdorff groupoid ${\mathcal G}$, and a locally constant $\Gamma $-valued cocycle $c$ on ${\mathcal G}$, let
  $$
   A_g = C_c({\mathcal G}_g,R), \quad \forall g\in \Gamma .
  $$
  Then the $ A_g$ are homogeneous submodules for a grading of $A_R(\mathcal G)$, which is moreover a symmetrical grading.
\end {prop}

\begin {proof}
  We leave it for the reader to verify that we indeed have a grading.  In order to prove symmetry, let $B$ be
  a compact open bissection in ${\mathcal G}_g$.  Then $B\inv \subseteq {\mathcal G}_{g\inv }$, so $1_{B\inv }\subseteq A_{g\inv }$.  Moreover, one has
  $$
  1_B = 1_B1_{B\inv }1_B\in A_g A_{g\inv } A_g.
  $$
  Since $ A_g$ is spanned by the elements of the form $1_B$, with $B$ as above, we conclude that
  $$
   A_g \subseteq A_g A_{g\inv } A_g,
  $$
  and, since the reverse inclusion is clearly true, the proof is complete.
\end {proof}

The fact that $A_R(\mathcal G)$ is symmetrically graded does not necessarily imply that all of its graded ideals share that
property.  However in the following special situation they do:

\begin {lemma}\label {equation:IdealsSymmGraded}
  Let $K$ be a
  field\/\footnote {From now on we will need our ring $R$ to be a field.  In this setting we write $K$ instead of $R$.},
  ${\mathcal G}$ be an ample, Hausdorff groupoid, and let $c$ be a locally constant $\Gamma $-valued cocycle on ${\mathcal G}$.  If $c\inv (e)$ is a
strongly effective groupoid, then every graded ideal of $A_K(\mathcal G)$ is symmetrically graded.  \end {lemma}

\begin {proof}
  Let $J$ be a graded ideal in $A_K(\mathcal G)$. Setting $J_g=J\cap A_g$, our task is to prove that
  $$
  J_g = J_gJ_{g\inv }J_g, \quad \forall g\in \Gamma .
  $$
  We leave the trivial inclusion
  ``$
  J_g \supseteq J_gJ_{g\inv }J_g
  $"
  for the reader to verify, and concentrate on the opposite one.
Given $g$ in $\Gamma $, pick any
  $
  f\in J_g,
  $
  and write
  $$
  f = \sum _{D \in F} a_D1_D,
  \eqno {(\dagger )}
  $$
  where $F$ is a finite collection of pairwise disjoint compact open bissections and the $a_D$ are nonzero elements of
$K$. 
By \cite[Lemma~2.2]{CE}, we may suppose, without loss of generally, that each $D$
in $F$ is contained in a single ${\mathcal G}_h$.

Denoting by $\pi _g$ the canonical projection from $A_K(\mathcal G)$ to $ A_g$, we then have
  $$
  f = \pi _g(f) =    \sum _{D \in F} a_D\pi _g(1_D).
  $$

  Notice that if $D$ is contained in ${\mathcal G}_h$, then $1_D$ lies in $ A_h$.  So
  $\pi _g(1_D)=1_D$, when $g=h$,
  while
  $\pi _g(1_D)=0$, otherwise.
  Consequently we have
  $$
  f = \sum _{\buildrel {\scriptstyle D \in F}\over {\vrule height 8pt width 0pt D\subseteq {\mathcal G}_g}} a_D1_D,
  $$
  which amounts to saying that, in the original description $(\dagger )$ of $f$, we may assume that $D\subseteq {\mathcal G}_g$, for every
$D$ in $F$.

  Our next goal will be to prove that $1_D$ lies in $J$, for each $D$ in $F$.
Given $D_0$ in $F$, set $f'=1_{D_0\inv }
f$.  Observing that $1_{D_0\inv }$ lies in $ A_{g\inv }$, we have
  $$
  f' = 1_{D_0\inv } f \in A_{g\inv } A_g \subseteq A_e,
  $$
  and it is clear that $f'$ also lies in $J$, so $f'\in J_e$.

Recalling that ${\mathcal G}_e $ is a strongly effective groupoid, we may
employ \cite [Lemma~4.3]{CEaHS} for the ideal
  $$
  J_e\trianglelefteq A_e = A_K({\mathcal G}_e),
  $$
  leading up to the conclusion that $f'|_{\go } \in J_{e}$.
  If $D$ is any member of $F$ other than $D_0$, hence disjoint from $D_0$, it is easy to see that $D_0\inv D$ does not
intercept $\go $, so
  $$
  f'|_{\go } = \Big (\sum _{D \in F} a_D1_{D_0\inv }1_D \Big )|_{\go } = \Big (\sum _{D \in F} a_D1_{D_0\inv D}\Big ) |_{\go } = a_{D_0}1_{D_0\inv D_0}.
  $$
  From this we deduce that
  $$
  1_{D_0} =
  (a_{D_0}\inv 1_{D_0})(a_{D_0} 1_{D_0\inv D_0}) =
  (a_{D_0}\inv 1_{D_0})(f'|_{\go }) \in J.
  $$
  This proves our claim that   $1_D$ lies in $J$, for every $D$ in $F$, and observe that $1_{D\inv }$ also lies in $J$
because
  $$
  1_{D\inv } =   1_{D\inv }   1_D   1_{D\inv } \in A_K(\mathcal G)J  A_K(\mathcal G)\subseteq J.
  $$

 Since $D\subseteq {\mathcal G}_g$, and hence $D\inv \subseteq {\mathcal G}_{g\inv }$, we have that $1_D\in J_g$, and
$1_{D\inv }\in J_{g\inv }$, so
  $$
  1_D = 1_D 1_{D\inv } 1_D \in J_g J_{g\inv } J_g,
  $$
  from where it follows that $f$ also lies in $J_g J_{g\inv } J_g$, thus proving that
  $$
  J_g\subseteq J_g J_{g\inv } J_g.
  \eqno \qedhere
  $$
\end {proof}

With this we may prove the following  main result:

\begin {thm}\label {thm:gradedideals}
  Let $K$ be a
  field,  ${\mathcal G}$ be an ample, Hausdorff groupoid, $\Gamma $ be a discrete group, and $c: {\mathcal G}\to \Gamma $ be a
locally constant cocycle such that
  $c^{-1}(e)$ is strongly effective.  Then the correspondence
  $$
  U \mapsto A_K({\mathcal G}|_U)
  $$
  is an isomophism from the lattice of open invariant subsets of $\go $ onto to the lattice of graded ideals in $A_K({\mathcal G})$.
\end {thm}

\begin {proof}
  Given any $f$ in $A_K(\mathcal G)$, and any $g$ in $\Gamma $, observe that the homogeneous  component $f_g$  
  of $f$ is obtained by taking the restriction of $f$ to the clopen set ${\mathcal G}_g$.
If $f$ is in $A_K({\mathcal G}|_U)$, where $U$ is a given open invariant subset of $\go $, it is then clear that $f_g$ is also
in $A_K({\mathcal G}|_U)$, so $A_K({\mathcal G}|_U)$ is indeed a graded ideal.  Therefore the correspondence mentioned in the
statement is well defined.

  Let $J$ be a graded ideal of $A_K(\mathcal G)$.  Then by
  \eqref {prop:GradedCorrespondence},
  \eqref {equation:SteinSymmGrded},
  \eqref {equation:SteinbSUnital} and
  \eqref {equation:IdealsSymmGraded},
  there is an invariant ideal
  $$
  I\trianglelefteq A_e = A_K\big ({\mathcal G}_e \big ),
  $$
  such that $J=  \text {Ind}(I)$.  Since ${\mathcal G}_e $ is assumed to be strongly effective, we have by
\cite [Theorem~3.1]{CEaHS}
that there is an open subset $U\subseteq \go $, which is invariant relative to ${\mathcal G}_e $, and such that $I$ consists of all $f$ in
$A_K\big ({\mathcal G}_e \big )$ whose support lies in
  $
  {\mathcal G}_e |_U.
  $
  In other words $I=A_K({\mathcal G}_e |_U)$.

We next claim that the invariance of $I$ (in the sense of \ref {def:Invariant}) implies that $U$ is in fact invariant
relative to ${\mathcal G}$ (in the sense of \ref {equation:InvariantSubset}).  In order to check this we will resort to \eqref
{equation:InvarCriteria}, so it suffices to prove that $A_K(\mathcal G)A_K(U)= A_K(U) A_K(\mathcal G)$.  Observing that
$\go \subseteq {\mathcal G}_e$, we have by
  \eqref {equation:NonGegenSubAlg} that
  $$
  A_K(\mathcal G)= A_K(\mathcal G)A_K({\mathcal G}_e) = A_K({\mathcal G}_e)   A_K(\mathcal G),
  $$
  so
  $$
  A_K(\mathcal G)A_K(U) =
  A_K(\mathcal G)A_K({\mathcal G}_e) A_K(U) \buildrel \eqref {equation:InvarCriteria} \over =
  A_K(\mathcal G)A_K({\mathcal G}_e|_U) =
  A_K(\mathcal G)I = $$$$ =
  I A_K(\mathcal G)=
  A_K({\mathcal G}_e|_U)  A_K(\mathcal G)=
  A_K(U) A_K({\mathcal G}_e)   A_K(\mathcal G)=
  A_K(U) A_K(\mathcal G),
  $$
  verifying the appropriate hypothesis in \eqref {equation:InvarCriteria}, and hence proving our claim that $U$ is invariant
relative to ${\mathcal G}$.

We next observe that
  $$
  J = \text {Ind}(I) = \medoplus _{g\in \Gamma }I A_g = I\big (\medoplus _{g\in \Gamma } A_g\big ) =
  IA_K(\mathcal G)=
  A_K({\mathcal G}_e|_U)  A_K(\mathcal G)=  $$$$ =
  A_K(U) A_K({\mathcal G}_e)  A_K(\mathcal G)=
  A_K(U)   A_K(\mathcal G)=
  A_K({\mathcal G}|_U).
  $$

This shows that our correspondence is surjective.
To see that it is also injective, let $U$ be an open invariant subset of $\go $ and observe that ${\mathcal G}|_U \cap \go =U$, so
  $$
  A_K({\mathcal G}|_U)\cap A_K(\go) = A_K(U).
  $$
  This says that $A_K(U)$, and hence also $U$, may be recovered from   $A_K({\mathcal G}|_U)$, from where
injectiveness follows.

  We leave it for the reader to prove that the corresponding lattice structures are preserved.
\end {proof}


 \section{Example:  Groups acting on graphs}

 \label{sec:gag}

In this section, we consider the algebras $\OGE$ associated to triples $(G,E, \varphi)$, introduced in \cite{EP_GGS}.

\subsection{The algebra $\OGE$}

Let us recall the construction.

\begin{noname}\label{basicdata}
{\rm The basic data for our construction is a triple $(G,E,\varphi)$ composed of:
\begin{enumerate}
\item A finite directed graph $E=(E^0, E^1, r, s)$ without sources.
\item A discrete group $G$ acting on $E$ by graph automorphisms.
\item A 1-cocycle $\varphi: G\times E^1\rightarrow G$ satisfying the property
$$
\varphi(g,a)\cdot v=g\cdot v \text{ for every } g\in G, a\in E^1, v\in E^0.
$$
\end{enumerate}
}
\end{noname}
The property $(3)$ required on $\varphi$ is tagged $(2.3)$ in \cite{EP_GGS}. 

\begin{definition}\label{Def:OGE}
{\rm Given a triple $(G,E, \varphi)$ as in  (\ref{basicdata}), we define $\OGE$ to be the universal $C^*$-algebra as follows:
\begin{enumerate}
\item \underline{Generators}:
$$\{p_x : x\in E^0\}\cup\{s_a : a\in E^1\} \cup \{u_g : g\in G\}.$$
\item \underline{Relations}:
\begin{enumerate}
\item $\{p_x : x\in E^0\}\cup\{s_a : a\in E^1\}$ is a Cuntz-Krieger $E$-family in the sense of \cite{Raeburn}.
\item The map $u:G\rightarrow \OGE$ defined by the rule $g\mapsto u_g$ is a unitary $\ast$-representation of $G$.
\item $u_gs_a=s_{g\cdot a}u_{\varphi(g,a)}$ for every $g\in G, a\in E^1$.
\item $u_gp_x=p_{g\cdot x}u_g$ for every $g\in G, x\in E^0$.
\end{enumerate}
\end{enumerate}
}
\end{definition}
\noindent Notice that the relation (2a) in Definition \ref{Def:OGE} implies that there is a natural representation map
$$
\begin{array}{cccc}
\phi: & C^*(E) &\to   & \OGE  \\
 & p_x & \mapsto  & p_x \\
 & s_a & \mapsto  & s_a  
\end{array}
$$
which is injective \cite[Proposition 11.1]{EP_GGS}.

\subsection{The groupoid $\CG$}

Recall from \cite[Definition 4.1]{EP_GGS} that given a triple $(G,E,\varphi)$ as in (\ref{basicdata}), we define a $\ast$-inverse semigroup $\SGE$ as follows:
\begin{enumerate}
\item The set is
$$\SGE=\{ (\alpha,g,\beta) : \alpha, \beta\in E^*, g\in G, d(\alpha)=gd(\beta)\}\cup \{ 0\},$$
where $E^*$ denotes the set of finite paths in $E$.
\item The operation is defined by:
$$(\alpha,g,\beta)\cdot (\gamma,h,\delta):=
\left\{
\begin{array}{cc}
 (\alpha, g\varphi(h,\varepsilon), \delta h\varepsilon), &   \text{if } \beta=\gamma \varepsilon  \\
   &  \\
 (\alpha g\varepsilon, \varphi(g,\varepsilon) h, \delta), &  \text{if } \gamma=\beta\varepsilon    \\
    &  \\
 0, & \text{otherwise,}
\end{array}
\right.
$$
and $(\alpha,g,\beta)^*:= (\beta,g^{-1}, \alpha)$.
\end{enumerate}

Then, we can construct the groupoid of germs of the action of $\SGE$ on the compact space $\widehat{\mathcal{E}}$
 of characters of the semilattice $\mathcal{E}$ of idempotents of $\SGE$. In our concrete case, $\widehat{\mathcal{E}}$ turns out to be 
 homeomorphic to the compact space $E^{\infty}$ of one-sided infinite paths on $E$; the action of 
 $(\alpha,g,\beta)\in \SGE$ on $\eta=\beta\widehat{\eta}$  is given by the rule $(\alpha,g,\beta)\cdot \eta=\alpha (g\widehat{\eta})$. 
Thus, the groupoid of germs is
$$\CG=\{[\alpha,g,\beta;\eta] : \eta=\beta\widehat{\eta}\},$$
where $[s;\eta]=[t;\mu]$ if and only if $\eta=\mu$ and there exists $0\ne e^2=e\in \SGE$ such that $e\cdot\eta=\eta$ and $se=te$.
The unit space 
\[\CG^{(0)}=\{[\alpha,1,\alpha;\eta] : \eta=\alpha\widehat{\eta}\}\]
is identified with the one-sided infinite path space $E^{\infty}$, via the
homeomorphism $[\alpha,1,\alpha; \eta]   \mapsto   \eta.$
Under this identification, the range and source maps on
$\CG$ are:
\[s([\alpha,g,\beta;\beta \widehat{\eta}])= \beta\widehat{\eta} \quad \text{ and } \quad  
r([\alpha,g,\beta;\beta \widehat{\eta}])= \alpha(g \widehat{\eta}).\]

A basis for the topology on $\CG$ is given by compact open bissections of the form 
\[ \Theta(\alpha,g,\beta;Z(\gamma)):= \{ [\alpha, g, \beta;\xi ] \in \CG : \xi \in  Z(\gamma) \} \] 
where $\gamma \in E^*$ and  $Z(\gamma):=\{\gamma \widehat{\eta} : \widehat{\eta}\in E^{\infty}\}.$
Thus $\CG$ is locally compact and ample.  
In \cite{EP_GGS} characterizations are given for when 
$\CG$ is Hausdorff \cite[Theorem 12.2]{EP_GGS}, 
amenable \cite[Corollary 10.18]{EP_GGS} or effective \cite[Theorem 14.10]{EP_GGS} 
in terms of the properties of the triple $(G,E,\varphi)$ and the action of $\SGE$ on $E^{\infty}$ 
(see \cite{ExPa} for a version of these results in the case of tight groupoids of 
arbitrary inverse semigroups). 

By \cite[Theorem 6.3 \& Corollary 6.4]{EP_GGS}, we have a $\ast$-isomorphism $\OGE\cong C^*(\CG)$, so that $\OGE$ can be seen as a full groupoid $C^*$-algebra. The complex Steinberg algebra $A_{\mathbb{C}}(\mathcal{G}_{(G,E)})$ is a dense subalgebra of $\mathcal{O}_{G,E} \cong C^*(\mathcal{G}_{(G,E)})$
by \cite[Proposition~6.7]{Steinberg}.

\subsection{The Steinberg algebra $A_R(\mathcal{G}_{(G,E)})$}

Now, we will  prove that for any unital commutative ring $R$, 
the Steinberg algebra $A_R(\mathcal{G}_{(G,E)})$ is isomorphic to the $R$-algebra 
$\mathcal{O}_{(G,E)}^{\text{alg}}(R)$ with presentation given by Definition \ref{Def:OGE} (i.e. the algebraic version of $\mathcal{O}_{(G,E)}$).

\begin{prop}\label{prop:TightRep}
The map
$$
\begin{array}{cccc}
\pi :& S_{G,E} & \to  & \mathcal{O}_{(G,E)}^{\text{alg}}(R)  \\
 & (\alpha, g, \beta)  & \mapsto  &   s_{\alpha}u_gs_{\beta}^*
\end{array}
$$
is the universal tight representation of $S_{G,E}$ in the category of $R$-algebras.
\end{prop}
\begin{proof}
This is a direct consequence of the arguments used to prove \cite[Proposition~6.2]{EP_GGS} and \cite[Theorem~6.3]{EP_GGS}.
\end{proof}

Since $\mathcal{G}_{(G,E)}$ is the groupoid of germs for the natural action of 
$S_{G,E}$ on the space of tight filters over its idempotent semi-lattice, we can 
argue as in the proof of \cite[Theorem~2.4]{E2} (see \cite[(2.4.2)]{E2}) to show that the map
$$\begin{array}{cccc}
\phi :& A_R(\mathcal{G}_{(G,E)}) & \rightarrow  & \mathcal{O}_{(G,E)}^{\text{alg}}(R)  \\
 & 1_{\Theta((\alpha, g, \beta), Z(\beta))}   & \mapsto  &  s_{\alpha}u_gs_{\beta}^*
\end{array}
$$
is a well-defined $R$-algebra homomorphism that is onto.

On the other side, the representation map
$$
\begin{array}{cccc}
\rho :& S_{G,E} & \rightarrow  & A_R(\mathcal{G}_{(G,E)})  \\
 & (\alpha, g, \beta)  & \mapsto  &   1_{\Theta((\alpha, g, \beta), Z(\beta))}
\end{array}
$$
is tight. Hence, the universal property stated in Proposition \ref{prop:TightRep} says
$$\begin{array}{cccc}
\psi :& \mathcal{O}_{(G,E)}^{\text{alg}}(R) & \rightarrow  & A_R(\mathcal{G}_{(G,E)})   \\
 &  s_{\alpha}u_gs_{\beta}^*  & \mapsto  &   1_{\Theta((\alpha, g, \beta), Z(\beta))}
\end{array}
$$
is an $R$-algebra homomorphism. Moreover, since $\widehat{\mathcal{E}(S_{G,E})}_{tight}$ is 
Boolean \cite{Steinberg}, $\psi$ is onto \cite[Proposition 5.13(7)]{Steinberg}. Clearly, $\phi$ and $\psi$ are mutually inverses. So, we conclude

\begin{thm}\label{thm:Iso1}
The map $\phi : A_R(\mathcal{G}_{(G,E)})  \rightarrow   \mathcal{O}_{(G,E)}^{\text{alg}}(R)$ is an $R$-algebra isomorphism. $\hspace{\fill} \Box$
\end{thm}

Notice that, through the isomorphism $\psi$, we can identify the Leavitt path algebra $L_R(E)$ of the graph $E$ 
(which is a subalgebra of $\mathcal{O}_{(G,E)}^{\text{alg}}(R)$) with a subalgebra of $A_R(\mathcal{G}_{(G,E)})$ 
isomorphic to the Steinberg algebra $A_R(\mathcal{G}_{E})$ of the path groupoid $\mathcal{G}_{E}$.\vspace{.2truecm}

An interesting application of Theorem \ref{thm:Iso1} is the following characterization of simplicity.
First we give the relavent definitions from \cite{EP_GGS}.
\begin{itemize}
\item We say $E$ is \emph{weakly $G$-transitive} if, given any infinite path $\zeta$, and any vertex $x \in E^0$, there is
some vertex $v$ along $\zeta$ such that there exists a vertex $u$ with $u=gx$ for some $g \in G$ 
and there is a path from $v$ to $u$.

\item A \emph{$G$-circuit} is a pair $(g, \gamma)$, where $g \in G$, and $\gamma \in E^*$ is a finite path of nonzero length
such that $s(\gamma) = gr(\gamma)$.

\item Given a $G$-circuit $(g,\gamma)$ such that
$\gamma = \gamma_1 \gamma_2 . . . \gamma_n$ in $E^*$ and each $\gamma_i$ is in $E^1$, 
we say that $\gamma$ has \emph{no
entry} if $r^{-1}(s(\gamma_i))$ is a singleton for every $i = 1, . . . , n$, i.e. $r^{-1}(s(\gamma_i))=\{\gamma_{i+1}\}$  for every $i = 1, . . . , n-1$, and $r^{-1}(s(\gamma_n))=\{g\gamma_{1}\}$. 

\item Given $g \in G$, and $x \in E^0$ , we shall say that $g$ is \emph{slack} at $x$, if there
is a non-negative integer $n$ such that all finite paths $\gamma$ with $r(\gamma) = x$, and $|\gamma| \geq n$, are
strongly fixed by $g$, as defined in the sense that $g\gamma = \gamma$, and $\phi(g, \gamma) = 1$.
\end{itemize}

\begin{cor}\label{cor:Simple1}
Suppose that $\mathcal{G}_{(G,E)}$ is Hausdorff. Then, for any field $K$ the following are equivalent:
\begin{enumerate}
\item The algebra $\mathcal{O}_{G,E}^{\text{alg}}(K)$ is simple.
\item The following properties hold:
\begin{enumerate}
\item $E$ is weakly $G$-transitive.
\item Every $G$-circuit has an entry.
\item Given $x\in E^0$ and $g\in G\setminus\{1\}$ fixing $Z(x)$ pointwise, then necessarily $g$ is slack at $x$.
\end{enumerate}
\end{enumerate}
\end{cor}
\begin{proof}
The result holds by Theorem \ref{thm:Iso1}, \cite[Theorem~4.1]{BCFS}, \cite[Theorem~13.6]{EP_GGS} and \cite[Theorem~14.10]{EP_GGS}.
\end{proof}

We are interested in determining the graded ideals of $A_K(\mathcal{G}_{(G,E)})$ for a suitable grading.
Notice that the function \[c:\mathcal{G}_{(G,E)} \to \Z\] such that $c([\alpha, g,\beta; \beta\xi]) = \vert\alpha\vert-\vert\beta\vert$
 is a groupoid cocycle and hence determines a $\Z$-grading on $A_{K}(\mathcal{G}_{(G,E)})$.
 However, in general $c^{-1}(0)$ will not be
 strongly effective so there might be $\Z$-graded ideals of $A_{K}(\mathcal{G}_{(G,E)})$ that
 have trivial intersection with $A_{K}(\go_{(G,E)})$.  In what follows, we show that 
 $A_{K}(\mathcal{G}_{(G,E)})$
 is also graded by the \emph{lag group}. 

 \subsection{The lag function}
Denote \[G^{\infty}:= \prod\limits_{n\in \mathbb{Z}^+} G \quad \text{ and } \quad
 G^{(\infty)}:= \bigoplus\limits_{n\in \mathbb{Z}^+} G.\]
Denote the \emph{Corona} group by
\[ \breve{G} := G^{\infty} / G^{(\infty)}. \]
Then the lag group of $\mathcal{G}_{(G,E)}$ is the semi-direct product $\group$ where $\breve{\rho}$
is defined as follows:
The \emph{right shift} is the endomorphism $\rho: G^{\infty} \to G^{\infty}$  such that
\[\quad \rho(\g)_n = \begin{cases}
                                                          \g_{n-1} & \text{if }n>0\\
                                                          1 &\text{if } n=0.
                                                         \end{cases}  \]
Since classes in $G^{\infty}$ are invariant under $\rho$,  we define $\breve{\rho}$ by passing to $\breve{G}$.
Since $\breve{\rho}$ is an automorphism on $\breve{G}$ (with an abuse of notation) we have
a homomorphism \[\breve{\rho}:\mathbb{Z} \to \operatorname{Aut}\breve{G}
\text{ such that } \breve{\rho}(n)(\g) = \breve{\rho}^n(\g).\]
\begin{notation}
We will write 1 for the identity in $G$, $\bf{1}$ for the identity in $\breve{G}$ and $\breve{\bf{1}}$ for the
identity in $\group$.
\end{notation}

Now, we define a map $\Phi: G\times E^{\infty}\to G^{\infty}$ by the rule 
$\Phi (g,\xi)_n:=\varphi(g, \xi\vert_{n-1})$ \cite[Definition~8.9]{EP_GGS}. Thus,

\begin{definition}[{\cite[Proposition~8.14]{EP_GGS}}]\label{def:Lag} 
The \emph{lag function} is the map
$$
\begin{array}{cccc}
\ell: & \mathcal{G}_{(G,E)} &  \to &  \group \\
 & [\alpha, g,\beta; \beta\xi] & \mapsto  &  (\breve{\rho}^{\vert \alpha\vert}(\breve{\Phi}(g, \xi)) ,\vert\alpha\vert-\vert\beta\vert)
\end{array}
$$
\end{definition}
The lag function is a one-cocycle by \cite[Proposition~8.14]{EP_GGS} and thus gives a grading on of our algebra.

\subsection{Graded ideals} We would like to apply Theorem \ref{thm:gradedideals} 
to determine the $\group$-graded ideals of $A_K(\mathcal{G}_{(G,E)})$ for any field $K$. To this end, we need to 
verify that $\ell^{-1}(e)$ is strongly effective.  In fact, we show it is principal.

\begin{lemma}\label{lem:FixedSubalgebra}
Suppose that $\mathcal{G}_{(G,E)}$ is Hausdorff. Let $e:=(\breve{\bf{1}},0)$ be the identity in $\group$. Then
$\ell^{-1}(e)=\{[{\alpha}, 1, {\beta}; \beta\xi] : \alpha, \beta\in E^*, \vert \alpha \vert =\vert \beta \vert\}$.
\end{lemma}
\begin{proof}
Let $\gamma\in \ell^{-1}(e)$. Then, $\gamma=[\alpha, g,\beta; \beta\xi]$, and since $\ell (\gamma)=(\breve{\bf{1}},0)$, we have that $\vert \alpha \vert =\vert \beta \vert$. On the other side, $\breve{\bf{1}}=\breve{\rho}^{\vert \alpha\vert}(\breve{\Phi}(g, \xi))$ means that, up to a finite number of entries, the sequence $(g, \varphi(g, \xi\vert _1), \varphi(g, \xi\vert _2), \dots ,\varphi(g, \xi\vert _k), \dots)$ coincides with the sequence $(1,1, \dots, 1, \dots)$. This is equivalent to say that there exists $n\in \mathbb{N}$ such that $\varphi(g, \xi\vert _k)=1$ for all $k\geq n$.  Define $\tau:=\xi\vert_n$. Then,
$$[\alpha, g,\beta; \beta\xi]=[\alpha\tau, \varphi(g, \tau),\beta\tau; \beta\tau\xi_{[n+1,\infty)}]=[\alpha\tau, 1,\beta\tau; \beta\tau\xi_{[n+1,\infty)}],$$
as desired.
\end{proof}

\begin{cor}\label{corol:NewLemma4.1}
Suppose that $\mathcal{G}_{(G,E)}$ is Hausdorff. Then, $\ell^{-1}(e)$ is principal.
\end{cor}
\begin{proof}
Clearly, $\go_{(G,E)}\subseteq \ell^{-1}(e)\cap \operatorname{Iso}(\mathcal{G}_{(G,E)})$. 
On the other side, if $\gamma \in \operatorname{Iso}(\mathcal{G}_{(G,E)})\setminus \go_{(G,E)}$, 
then $\gamma=[\alpha, g,\alpha, \alpha\xi]$ with $g\cdot \xi=\xi$ but $\varphi(g,\xi\vert_n)\ne 1$ 
for every $n\in \mathbb{N}$. By Lemma \ref{lem:FixedSubalgebra}, $\gamma \not\in \ell^{-1}(e)$, so we are done.
\end{proof}

Thus we can apply Theorem~\ref{thm:gradedideals} to see that the graded ideals
in $A_{\C}(\mathcal{G}_{(G,E)})$ are precisely the ideals of the form let $I(U):= A_{\C}(\mathcal{G}_{(G,E)|U})$
for some open $\SGE$-invariant subset $U \subseteq \go_{(G,E)}$.

Next we will show how graded ideals in $A_{\C}(\mathcal{G}_{(G,E)})$ 
come from particular subsets of vertices of $E$.
We say $H \subseteq E^0$ is \emph{hereditary} if for any $w \in H$ such that there exists $e \in E^1$
with $r(e) = w$,  then $s(e) \in H$ as well.  We say $H$ is \emph{saturated} if the following condition is satisfied:
\[v \in E^0 \text{ such that } s(r^{-1}(v)) \subseteq H \implies v \in H.\]
Now suppose that $U$ is an open $\SGE$-invariant subset of $\go_{G,E} \cong E^{\infty}$.
Motiated by the graph algebra construction in \cite[Theorem~3.3]{CMMS} (for example) define
\[
 H_U := \{x \in E^0 : \Theta((x,1,x), Z(x)) \subseteq U\}.
\]
\begin{lemma}
 The set $H_U$ defined above is a $G$-invariant saturated hereditary set.
\end{lemma}

\begin{proof}
To see that $H_U$ is hereditary, fix $w \in H_U$ and $e \in E^1$ such that
$r(e) = w$.  It suffices to show that $\Theta(s(e),1,s(e); Z(s(e))) \subseteq U$.
Fix \[[s(e), 1, s(e) ; \zeta] \in \Theta(s(e),1,s(e); Z(s(e))).\]  Notice 
\begin{equation}
\label{eqn:HU}
s([e, 1, s(e) ; \zeta]) = [s(e), 1, s(e) ; \zeta] \quad \text{ and } \quad r([e, 1, s(e) ; \zeta]) = [e, 1, e ; e\zeta].
\end{equation}
Since $[e, 1, e ; e\zeta] \in \Theta(w,1,w ; Z(w)) $ we have $[e, 1, e ; e\zeta] \in U$ and hence 
$[s(e), 1, s(e) ; \zeta] \in U$ because $U$ is invariant.

To see that $H_U$ is saturated, fix $v$ such that $s(r^{-1}(v)) \subseteq H_U$. 
Notice that \[\Theta(v,1,v; Z(v)) = \bigcup_{e \in r^{-1}(v)} \Theta(e,1,e; Z(e)). 
            \]
            We have $\Theta(s(e),1,s(e); Z(s(e))) \subseteq U$.  Using an argument similarly to the one in \eqref{eqn:HU},
            we see that $\Theta(e,1,e; Z(e)) \subseteq U$.
            Thus $\Theta((v,1,v), Z(v)) \subseteq U$ which means $v \in H_U$.
            
That $H_U$ is $G$-invariant follows from $U$ being $\SGE$-invariant.             

\end{proof}

Denote by $1_x$ the characteristic function $1_{\Theta((x,1,x), Z(x))}\in A_{\C}(\CG)$. 
Also, let $I(H_U)$ denote the ideal generated by $\{1_x\}_{x \in H_U}$. 
\begin{prop}
\label{Prop:D<->H_D}
Suppose that $\mathcal{G}_{(G,E)}$ is Hausdorff. 
Let $U$ be an open $\SGE$-invariant subset in $E^{\infty}$.  Then
 $I(U) = I(H_U)$. 
\end{prop}

\begin{proof}
To  show the forward containment, it suffices to show $1_{\Theta((\alpha,1,\alpha), Z(\alpha))}$ is in $I(H_U)$ 
for any $\Theta((\alpha,1,\alpha), Z(\alpha)) \subseteq U$.
First we claim that  $\Theta((s(\alpha),1,s(\alpha)), Z(s(\alpha))) \subseteq U$.
To see this, fix \[[s(\alpha),1,s(\alpha); \zeta] \in \Theta((s(\alpha),1,s(\alpha)), Z(s(\alpha)).\]
Then we have
\[
s([\alpha,1,s(\alpha); \zeta]) = [s(\alpha),1,s(\alpha); \zeta] \quad \text{ and } 
\quad r([\alpha,1,s(\alpha); \zeta]) = [\alpha,1,\alpha; \alpha \zeta] \in U.
\]
Thus $[s(\alpha),1,s(\alpha); \zeta] \in U$ by invariance, proving the claim.
Therefore $s(\alpha) \in H_U$.  Now \[1_{\Theta((\alpha,1,\alpha), Z(\alpha))} = 1_{\Theta((\alpha,1,\alpha), Z(\alpha))} 1_{s(\alpha)}\] 
and hence $1_{\Theta((\alpha,1,\alpha), Z(\alpha))} \in I(H_U).$ 

For the reverse containment, notice that for any $x \in H_U$ we have $1_x \in I(U)$. 
\end{proof}

As an immediate consequence of Proposition \ref{Prop:D<->H_D}, we conclude that
\begin{cor}\label{RenIdealsAreGaugeInv}
If $U$ is an open $\SGE$-invariant subset of $E^{\infty}$, then $I(U)$ is a graded ideal of 
$A_K(\CG)$ generated by a $G$-invariant saturated hereditary subset 
$H_U$ of $E^0$. \hspace{\fill} $\square$
\end{cor}

Thus, the relevant information related to this concrete class of ideals relates to the ideal theory developed 
for graph $C^*$-algebras \cite[Chapter 4]{Raeburn} and Leavitt path algebras \cite[Proposition 5.2 \& Theorem 5.3]{AMFP}.

%
%
%

 \section{Example:  Boolean dynamical systems}
 \label{sec:bds}

$C^*$-algebras associated to Boolean dynamical systems were recently introduced by Carlsen, Ortega and the third author in \cite{COP}.
Let us briefly recall the definition of these algebras.

\subsection{The algebra}
Let $\B$ be a Boolean algebra, we say that a map $\theta:\B\longrightarrow \B$ is an \emph{action on $\B$} 
if $\theta$ is a Boolean algebra homomorphism with $\theta(\emptyset)=\emptyset$.
We say that the action has \emph{compact range} if $\{\theta(A)\}_{A\in \B}$ has least upper-bound, 
that we will denote $\mathcal{R}_\theta$. Moreover, we say that the action has \emph{closed domain}  
if there exists $\mathcal{D}_\theta\in \B$ such that $\theta(\mathcal{D}_\theta)=\mathcal{R}_\theta$.

Given a set $\Labe$, and given any $n\in\N$, we define 
\[\Labe^n=\{(\alpha_1,\ldots,\alpha_n):\alpha_i\in \Labe)\} \quad \text{ and } \quad \Labe^*=\bigcup\limits_{n=0}^\infty\Labe^n,\] where $\Labe^0=\{\emptyset\}$.  Given $\alpha\in \Labe^n$ for $n\geq 1$, 
we will write it as $\alpha=\alpha_1\cdots\alpha_n$ where $\alpha_i\in \Labe$. \vspace{.2truecm}

A \emph{Boolean dynamical system} on a Boolean algebra $\B$ is a triple $(\B,\Labe,\theta)$ such that $\Labe$ is a set, 
and $\{\theta_\alpha\}_{\alpha\in \Labe}$ is a set of actions on $\B$. Moreover, given $\alpha=(\alpha_1,\ldots,\alpha_n)\in \Labe^{\geq 1}$ 
the action $\theta_\alpha:\B\longrightarrow \B$ defined as $\theta_\alpha=\theta_{\alpha_n}\circ\cdots\circ\theta_{\alpha_1}$ has compact range and closed domain.

Given any $\alpha\in \Labe^*$, we will write $\mathcal{D}_\alpha:=\mathcal{D}_{\theta_\alpha}$ and 
$\mathcal{R}_\alpha:=\mathcal{R}_{\theta_\alpha}$. Also, when $\alpha=\emptyset$, we will define 
$\theta_{\emptyset}=\text{Id}$, and we will formally assume that $\mathcal{R}_{\emptyset}=\mathcal{D}_{\emptyset}:=\bigcup\limits_{A\in \mathcal{B}}A$, in order to guarantee that $A\subseteq \mathcal{R}_{\emptyset}$ for every $A\in \mathcal{B}$.

Let $(\B,\Labe,\theta)$ be a  Boolean dynamical system. Given $B\in \B$ we define
$$\Delta_{B}:=\{\alpha\in \Labe: \theta_\alpha(B)\neq\emptyset\}\qquad\text{and}\qquad \lambda_{B}:=|\Delta_{B}|\,.$$  
We say that $A\in \B$ is a \emph{regular set} if given any $\emptyset\neq B\in \B$ with 
$B\subseteq A$ we have that $0<\lambda_B<\infty $, otherwise $A$ is called a \emph{singular set}. 
We denote by $\Breg$ the set of all regular sets where we will include $\emptyset$.\vspace{.2truecm}

A \emph{Cuntz-Krieger representation of the Boolean dynamical system} $(\B,\Labe, \theta)$ in a 
$C^*$-algebra $\mathcal{A}$ consists of a family of projections $\{P_A : A\in\B\}$ and partial isometries 
$\{S_\alpha : \alpha\in \Labe\}$ in $\mathcal{A}$, with the following properties:
\begin{enumerate}
\item If $A,B\in\B$, then $P_A\cdot P_B = P_{A\cap B}$ and $P_{A\cup B} = P_A + P_B-P_{A\cap B}$, where $P_{\emptyset}=0$.
\item If $\alpha\in \Labe$ and $A\in\B$, then $P_A\cdot  S_\alpha = S_\alpha  \cdot  P_{\theta_\alpha(A)}$.
\item If $\alpha, \beta\in \Labe$ then $S_\alpha^* \cdot S_\beta = \delta_{\alpha,\beta} \cdot  P_{\Ra}$.
\item Given $A\in \Breg$ we have that
$$P_A =\sum_{\alpha\in\Delta_A} S_\alpha\cdot  P_{\theta_\alpha(A)}\cdot  S_\alpha^*\,.$$
\end{enumerate}
A representation is called \emph{faithful} if $P_A\neq 0$ for every $A\in\B$.

Given a representation $\{P_A,S_\alpha\}$ of a Boolean dynamical system $(\B,\Labe, \theta)$ in 
a $C^*$-algebra $\mathcal{A}$, we define $C^*(P_A,S_\alpha)$ to be the sub-$C^*$-algebra
of $\mathcal{A}$ generated by $\{P_A,S_\alpha:A\in \B,\,\alpha\in \Labe\}$.

A \emph{universal representation} $\{p_A,s_\alpha\}$ of a Boolean dynamical system $(\B,\Labe, \theta)$ is a
representation satisfying the following universal property: given a representation $\{P_A,S_\alpha\}$ of
$(\B,\Labe,\theta)$ in a $C^*$-algebra $\mathcal{A}$, there exists a non-degenerate $*$-homomorphism
$\pi_{S,P}:C^*(p_A,s_\alpha)\longrightarrow \mathcal{A}$ such that $\pi_{S,P}(p_A)=P_A$ and $\pi_{S,P}(s_\alpha)=S_\alpha$ for
$A\in \B$ and $\alpha\in \Labe$.  We will set $C^*(\B,\Labe,\theta):=C^*(p_A,s_\alpha)$.\vspace{.2truecm}

It is shown that
$$T=T_{(\B,\Labe,\theta)}:=\{s_\alpha p_A s^*_\beta: \alpha,\beta\in \Labe^*\,, A\in \B\,, A\subseteq \Ri_\alpha\cap\Ri_\beta\neq 
\emptyset\}\cup\{0\}\subseteq \Clab\,,$$
is a $\ast$-inverse semigroup \cite[Proposition~6.2]{COP} such that $\Clab$ is the closure of the linear span of
$T$, and that it is $E^*$-unitary \cite[Proposition~6.8]{COP}. Moreover, if $\B$ and $\Labe$ are countable sets, then 
$\Clab$ is $\ast$-isomorphic to $C^*(\mathcal{G}_{tight}(T))$ \cite[Theorem~8.3]{COP}, which is a Hausdorff 
\cite[Lemma~8.2]{COP} amenable \cite[Lemma~8.4]{COP} ample groupoid. 
\vspace{.3truecm}

\subsection{The Steinberg algebra}
An argument analogous to that used for proving Theorem \ref{thm:Iso1} shows

\begin{thm}\label{thm: iso2}
The map
$$
\begin{array}{cccc}
\tau : & A_R(\mathcal{G}_{tight}(T))  &  \rightarrow &  R(\B, \Labe, \theta) \\
 & 1_{\Theta(s_{\alpha}, \operatorname{Dom}(s_{\alpha}^*s_{\alpha}))}  &  \mapsto &   s_{\alpha}\\
  & 1_{\Theta(p_A, \operatorname{Dom}(p_A))}  &  \mapsto &   p_A
\end{array}
$$
is an $R$-algebra isomorphism. $\hspace{\fill} 	\Box$
\end{thm}

In particular, we can characterize simplicity of $A_R(\mathcal{G}_{tight}(T))$. To do this, let us recall some definitions

Let $(\B,\Labe,\theta)$ be a Boolean dynamical system. Then:
\begin{itemize}
\item We say that the pair $(\alpha,A)$ with  $\alpha=\alpha_1\cdots\alpha_n\in \Labe^n$ for some $n\geq 1$, 
and $\emptyset\neq A\in \B$ with $A\subseteq\Ri_\alpha$, is a \emph{cycle} 
if given $k\in\N\cup\{0\}$ we have that $\theta_{\alpha^k}(A)\neq \emptyset$ and for every 
$\emptyset\neq B\subseteq \theta_{\alpha^k}(A)$ we have that $B\cap \theta_{\alpha}(B)\neq \emptyset$.
\item A cycle $(\alpha,A)$ has \emph{no exits} if  
given any $k\in\N\cup\{0\}$ we have that $\theta_{\alpha^k\alpha_1\cdots\alpha_t}(A)\in\Breg$ with  
$\Delta_{\theta_{\alpha^k\alpha_1\cdots\alpha_t}(A)}=\{\alpha_{t+1}\}$ for $t<n$ and $\theta_{\alpha^{k+1}}(A)\in\Breg$ with $\Delta_{\theta_{\alpha^{k+1}}(A)}=\{\alpha_1\}$.
\item We say that $(\B,\Labe,\theta)$ satisfies \emph{condition $(L_\B)$} if there is no cycle without exits.
\item We say that an ideal $\mathcal{I}$ of $\B$ is \emph{hereditary} if given $A\in \mathcal{I}$ and $\alpha\in \Labe$ then 
$\theta_\alpha(A)\in \mathcal{I}$. 
\item We also say that $\mathcal{I}$ is \emph{saturated} if given $A\in \Breg$  with 
$\theta_{\alpha}(A)\in \mathcal{I}$ for every $\alpha\in \Delta_A$ then $A\in \mathcal{I}$.
\end{itemize}
\vspace{.2truecm}

As a consequence, we have the following:

\begin{cor}\label{cor:SimpleBoolean}
If $(\B, \Labe, \theta)$ is a Boolean dynamical system such that $\B$ and 
$\mathcal{L}$ are countable, $K$ is a field and $K(\B,\Labe, \theta)$ is its associated $K$-algebra, then, the following statements are equivalent:
\begin{enumerate}
\item $K(\B,\Labe, \theta)$ is simple.
\item The following properties hold:
\begin{enumerate}
\item $(\B,\Labe,\theta)$ satisfies condition $(L_\B)$, and
\item The only hereditary and saturated ideals of $\B$ are $\emptyset$ and $\B$.
\end{enumerate}
\end{enumerate}
\end{cor}
\begin{proof}
The result holds by Theorem \ref{thm: iso2}, \cite[Theorem~4.1]{BCFS}, \cite[Theorem~9.7]{COP} and \cite[Theorem~9.15]{COP}.
\end{proof}

\vspace{.3truecm}

\subsection{Graded ideals}
For any field $K$, we can endow a structure of $\Z$-graded algebra on $A_K(\mathcal{G}_{\text{tight}}(T))$ using the continuous cocycle
$$
\begin{array}{cccc}
c:  & \mathcal{G}_{tight}(T) & \rightarrow   & \mathbb{Z}  \\
 & [s_{\alpha}p_As_{\beta}^*, \xi] & \mapsto  & \vert\alpha\vert -\vert\beta\vert
\end{array}.
$$
Hence, we apply Theorem \ref{thm:gradedideals} to determine the $\Z$-graded
ideals of $A_K(\mathcal{G}_{\text{tight}}(T))$ for any field $K$, provided we are able to prove that $c^{-1}(0)$ is strongly effective. We will prove that, in fact, $c^{-1}(0)$ is  principal.

\begin{lemma}\label{lem:BooleanPrincipal}
The subgroupoid $c^{-1}(0)$ of $\mathcal{G}_{tight}(T)$ is principal.
\end{lemma}
\begin{proof}
First, we define the $\ast$-subsemigroup of $T$
$$S=\{s_{\alpha}p_As_{\beta}^*\in T : \vert\alpha\vert =\vert\beta\vert\}.$$
Miming the arguments in \cite[Section 7]{COP}, one can show that
$$\iota : S\rightarrow ({\Clab^{\gamma}})^1,$$
is the universal tight representation of $S$ (where $\gamma: \mathbb{T}\rightarrow \Clab$ is the standard gauge action),
and that $\Clab^{\gamma}\cong C^*(\mathcal{G}_{tight}(S))$. Moreover, it is a simple exercise to
prove that $\mathcal{G}_{tight}(S)$ and $c^{-1}(0)$ are topologically isomorphic.

Then, an adaptation of the proof of \cite[Proposition 2.9]{OrtegaPriv} 
(see also the proof of \cite[Theorem 4.4]{BP}or \cite[Lemma 2.2]{BPRS}) shows that $\Clab^{\gamma}$ is a AF-algebra; 
we thank Eduard Ortega for turning our attention to this fact, and providing an accurate proof of it. 
Thus, \cite[Proposition III.1.15 \& Remarks III.1.2]{Ren} imply that  $c^{-1}(0)$ is principal, as desired.
\end{proof}

We then apply Theorem \ref{thm:gradedideals} to determine the $\Z$-graded
ideals of $A_K(\mathcal{G}_{\text{tight}}(T))$ for any field $K$. First, we recall some 
more definitions from \cite{COP}.
Given a collection $\Id$ of elements of $\B$ we define the \emph{hereditary expansion} of $\Id$ as
$$\He(\Id):=\{B\in \B: B\subseteq \bigcup\limits_{i=1}^n \theta_{\alpha_i}(A_i)\text{ where }A_i\in \Id \text{ and }\alpha_i\in \Labe^*\}\,.$$
Clearly, $\He(\Id)$ is  the minimal hereditary ideal of $\B$ containing $\Id$. Also, we define the \emph{saturation} of $\Id$, 
denoted by $\Sa(\Id)$, to be the minimal ideal of $\B$ generated by the set
$$\bigcup\limits_{n=0}^\infty\Sa^{[n]}(\Id),$$
defined by recurrence on $n\in\Z^+$ as follows:
\begin{enumerate}
\item $\Sa^{[0]}(\Id):=\Id$.
\item For every $n\in \N$, $\Sa^{[n]}(\Id):=\{B\in \Breg:\theta_\alpha(B)\in\Sa^{[n-1]}(\Id)\text{ for every }\alpha\in \Delta_B\}\,$.
\end{enumerate}
Observe that if $\Id$ is hereditary, 
then $\Sa(\Id)$ is also hereditary. Therefore, given a collection $\Id$ 
of elements of $\B$, $\Sa(\He(\Id))$ is the minimal hereditary and saturated ideal of $\B$ containing $\Id$.

Given any idempotent $e\in \mathcal{E}(T)$, let $D_e$ denote the domain of the (partial) 
action of $e$ on the space of tight filters of $\mathcal{E}(T)$, which is homeomorphic to $\mathcal{G}_{\text{tight}}(T)^{(0)}$. 
Also, we denote by $\mathcal{O}(e):=r(s^{-1}(D_e))$ the orbit of $D_e$ by the action of $T$, which is an 
open invariant subset of $\mathcal{G}_{\text{tight}}(T)^{(0)}$. Since $\{D_e : e\in \mathcal{E}(T)\}$ is a 
basis of compact open sets of $\mathcal{G}_{\text{tight}}(T)^{(0)}$, then $\{\mathcal{O}(e) : e\in \mathcal{E}(T)\}$ 
is a cover of the collection of open invariant sets of $\mathcal{G}_{\text{tight}}(T)^{(0)}$ such that 
any open invariant subset $U$ of $\mathcal{G}_{\text{tight}}(T)^{(0)}$ is a union of sets in this collection. 
Notice that, given any $e\in \mathcal{E}(T)$, $e$ is of the form $s_{\alpha}p_As_{\alpha}^*$ for 
$\alpha \in \mathcal{L}^*$ and $A\subseteq \mathcal{R}_{\alpha}$. Thus, taking $s=s_{\alpha}^*$, 
we have that $s\cdot(D_e)=D_{ses^*}$, where $ses^*=p_A$. Hence, for any $e\in \mathcal{E}(T)$ there 
exists a (unique) $A\in \mathcal{B}$ such that $\mathcal{O}(e)=\mathcal{O}(p_A)$.\vspace{.2truecm}

In an analogous way to that of \cite[Section 10]{COP}, we introduce the following definitions.

\begin{definition}\label{Def:Grad-HerSat}
Let $(\B, \Labe, \theta)$ be a (countable)  Boolean dynamical system. Then:
\begin{enumerate}
\item Given $U\subseteq \mathcal{G}_{\text{tight}}(T)^{(0)}$ an open invariant subset, we define
$$\He_U:=\{A\in \mathcal{B} : D_{p_A}\in I(U)\}$$
where $I(U)$ is the ideal of $A_K(\mathcal{G}_{\text{tight}}(T))$ generated by $U$. Clearly, $\He_U$ is a hereditary and saturated subset of $\mathcal{B}$.
\item Given any hereditary and saturated subset $\He$ of $\mathcal{B}$, we define
$$U_{\He}:=\bigcup\limits_{A\in \He}\mathcal{O}(p_A)$$
which is clearly an open invariant subset of $\mathcal{G}_{\text{tight}}(T)^{(0)}$.
\end{enumerate}
\end{definition}

\begin{prop}\label{Prop:Bijection}
Let $(\B, \Labe, \theta)$ be a (countable)  Boolean dynamical system. Then:
\begin{enumerate}
\item Given $U\subseteq \mathcal{G}_{\text{tight}}(T)^{(0)}$ an open invariant subset, $U=U_{\He_U}$.
\item Given any hereditary and saturated subset $\He$ of $\mathcal{B}$, $\He=\He_{U_\He}$.
\end{enumerate}
So, there is a lattice isomorphism between open invariant subsets of $\mathcal{G}_{\text{tight}}(T)^{(0)}$ and hereditary and 
saturated subset of $(\B, \Labe, \theta)$.
\end{prop}
\begin{proof} $\mbox{ }$
(1) Take $U\subseteq \mathcal{G}_{\text{tight}}(T)^{(0)}$ an open invariant subset. Since $A\in \He_U$ if and only if $D_{p_A}\subseteq U$ and $U$ is invariant,
we have that $\mathcal{O}(p_A)\subseteq U$. Thus, as $U_{\He_U}=\bigcup\limits_{A\in \He_U}\mathcal{O}(p_A)$, we have that $U_{\He_U}\subseteq U$. 
On the other side, if $e\in T$ is an idempotent and $D_e\subseteq U$, then there exists a unique $A\in \mathcal{B}$ such 
that $\mathcal{O}(e)=\mathcal{O}(p_A)$, whence $D_{p_A}\subseteq U$. Thus, $D_e\subseteq \mathcal{O}(e)=\mathcal{O}(p_A)\subseteq U_{\He_U}$. 
Hence, $U\subseteq U_{\He_U}$, so we are done.

(2) Given any hereditary and saturated subset $\He$, we have that
$$\He_{U_\He}=\{A\in \mathcal{B} : D_{p_A}\subseteq \bigcup\limits_{B\in \He}\mathcal{O}(p_B)\}.$$
Thus, $\He \subseteq \He_{U_\He}$ trivially. On the other side, if $A\in \He_{U_\He}$, then $D_{p_A}\subseteq \bigcup\limits_{B\in \He}\mathcal{O}(p_B)$. 
Since $D_{p_A}$ is compact, there exists $B_1, \dots, B_n\in \He$ such that $D_{p_A}\subseteq \bigcup\limits_{i=1}^{n}\mathcal{O}(p_{B_i})$. Since $\mathcal{O}(p_B)=\bigcup\limits_{s\in T}D_{sp_Bs^*}$, again by compactness, for each $1\leq i\leq n$ there exists $s_{i1}, \dots, s_{i{m_i}}\in T$ such that $D_{p_A}\subseteq \bigcup\limits_{i=1}^{n}\bigcup\limits_{j=1}^{m_i}D_{s_{ij}p_{B_i}s_{ij}^*}$. Since the elements of $T$ are of the form $s_{\alpha}p_As_{\beta}^*$ with $\alpha, \beta \in\mathcal{L}^*$ and $B\subseteq \mathcal{R}_{\alpha}\cap\mathcal{R}_{\beta}$, we can assume, relabelling if necessary, that there exist $\bigcup\limits_{i=1}^n\{\alpha_{i1},\dots ,\alpha_{i{m_i}}\}\subseteq \mathcal{L}^*$ such that
$$D_{p_A}\subseteq \bigcup\limits_{i=1}^{n}\bigcup\limits_{j=1}^{m_i}D_{s_{\alpha_{ij}}p_{\theta_{\alpha_{ij}}(B_i)}s_{\alpha_{ij}}^*}.$$
By \cite[Proposition 3.7]{ExPa}, this is equivalent to say that $\bigcup\limits_{i=1}^n\bigcup\limits_{j=1}^{m_i}\{s_{\alpha_{ij}}p_{\theta_{\alpha_{ij}}(B_i)}s_{\alpha_{ij}}^*\}$ is an outer cover of $p_A$. 
Then, applying the same argument as in \cite[Theorem  9.5, proof of (3)$\Rightarrow$ (1)]{COP}, we conclude that there exists $N\in\N$ such that $A\in \Sa^{[N]}(\He (B_1, \dots,B_N))\subseteq \He$. Hence, $\He_{U_\He}\subseteq \He$, so we are done.
\end{proof}

As a consequence of Proposition \ref{Prop:Bijection} and Theorem \ref{thm:gradedideals}, we conclude

\begin{thm}\label{Thm:IdealOfBoolDynSys}
Given a (countable) Boolean dynamical system $(\B, \Labe, \theta)$ and a field $K$, 
the lattice of $\Z$-graded ideals of $A_K(\mathcal{G}_{\text{tight}}(T))$ is isomorphic to the lattice of hereditary and saturated subset of $\Clab$.
\end{thm}

Observe that Theorem \ref{Thm:IdealOfBoolDynSys} is analog to \cite[Theorem 10.12]{COP}, where the result is proved for $C^*$-algebras of Boolean dynamical systems.



\begin{thebibliography}{99}

\bibitem{AMFP} \textsc{P. Ara, M.A. Moreno, E. Pardo}, Nonstable K-Theory for graph algebras, \emph{Algebras Represent. Theory.} \textbf{10} (2007), 157--178.

\bibitem{APS} \textsc{G. Aranda-Pino, E. Pardo, M. Siles-Molina}, Exchange Leavitt path algebras and stable rank, \emph{J. Algebra} \textbf{305} (2006), 912--936.

\bibitem{BP} \textsc{T. Bates, D. Pask}, $C^*$-algebras of labelled graphs. II. Simplicity results, \emph{Math. Scand.} \textbf{104} (2009), no. 2, 249--274.

\bibitem{BPRS} \textsc{T. Bates, D. Pask, I. Raeburn, W. Szyma\'{n}ski}, The $C^*$-algebras of row-finite graphs, \emph{New York J. Math.} \textbf{6} (2000), 307--324.

\bibitem{BCFS} \textsc{J. Brown, L.O. Clark, C. Farthing, A. Sims}, Simplicity of algebras associated to \'etale groupoids,
\emph{Semigroup Forum} \textbf{88} (2014), 433--452.

\bibitem{BNRSW} \textsc{J.H. Brown, G. Nagy, S. Reznikoff, A. Sims and D.P. Williams},
Cartan subalgebras in $C^*$-algebras of Hausdorff \'etale groupoids,  \emph{Integral Equations Operator Theory} \textbf{85} (2016), 109--126.

\bibitem{COP} \textsc{T.M. Carlsen, E. Ortega and E. Pardo}, $C^*$-algebras associated to Boolean dynamical systems, arXiv:1510.06718v2.

\bibitem{CE} \textsc{L.O. Clark and C. Edie-Michell}, Uniqueness theorems for Steinberg algebras, \emph{Algebr. Represnt. Theor.} \textbf{18} (2015) 907--916.

\bibitem{CEaHS} \textsc{L.O. Clark, C. Edie-Michell, A. an Huef and A. Sims}, Ideals of Steinberg algebras of strongly effective groupoids,
with applications to Leavitt path algebras, arXiv:1601.07238.


\bibitem{CFST} L.O. Clark, C. Farthing, A. Sims and M. Tomforde, A groupoid generalization of Leavitt path algebras,
\emph{Semigroup Forum}, \textbf{89} (2014), 501--517.

\bibitem{CGN} L.O. Clark, C. Gil Canto and A. Nasr-Isfahani, The cycline subalgebra of a Kumjian-Pask algebra,
 arXiv:1603.00508.

\bibitem{CMMS} L.O. Clark, D. Mart\'in Barquero, C. Mart\'in Gonz\'alez and M. Siles Molina,
Using Steinberg algebras to study decomposability of Leavitt path algebras, arXiv;1603:1033v1.

\bibitem{E2011} \textsc{R.Exel}, Non-Hausdorff \'etale groupoids, \emph{Proc. Amer. Math. Soc.} \textbf{139} (2011), 897--907.

\bibitem{E} \textsc{R. Exel}, Inverse semigroups and combinatorial $C^*$-algebras,
\emph{Bull. Braz. Math. Soc.} \textbf{39} (2008),
191--313.

\bibitem{E2} \textsc{R. Exel}, Reconstructing a totally disconnected groupoid from its ample semigroup,
\emph{Proc. Amer. Math. Soc.} \textbf{138} (2008), 2991--3001.

\bibitem{ExPa} \textsc{R. Exel, E. Pardo}, The tight groupoid of an inverse semigroup,
\emph{Semigroup Forum} \textbf{92} (2016), 274--303.

\bibitem{EP_GGS} \textsc{R. Exel and E. Pardo}, Self-similar graphs, a unified treatment
of Katsura and Nekrashevych $C^*$-algebras, arXiv:1409.1107v2.

\bibitem{GN} \textsc{C. Gil Canto, A. Nasr-Isfahani}, The maximal commutative subalgebra of a
Leavitt path algebra, arXiv: 1510.03992.


\bibitem{OrtegaPriv} \textsc{E. Ortega}, Simple Cuntz-Krieger Boolean algebras, \emph{Draft} (2016).

\bibitem{Raeburn} \textsc{I. Raeburn}, ``Graph Algebras'', CBMS Reg. Conf. Ser. Math., vol. 103, Amer. Math.
Soc., Providence, RI, 2005.

\bibitem{Ren} \textsc{J. Renault}, ``A groupoid approach to $C^*$-algebras'', \emph{Lecture Notes in Mathematics} \textbf{793}, Springer, Berlin, 1980.

\bibitem{Steinberg} \textsc{B. Steinberg}, A groupoid approach to discrete inverse semigroup algebras,
\emph{Adv. Math.} \textbf{223} (2010) 689--727.

\bibitem{St2016} \textsc{B. Steinberg},  Simplicity, primitivity and semiprimitivity of \'etale groupoid algebras with
applications to inverse semigroup algebras,  \emph{J. Pure Appl. Algebra} \textbf{220} (2016), 1035--1054.
\end{thebibliography}
\end{document}